\documentclass[a4paper,10pt,reqno]{amsart}
\usepackage{amsmath,amsthm,amssymb}
\usepackage{graphicx, color}
\usepackage[numbers,sort&compress]{natbib}
\nonstopmode \numberwithin{equation}{section}
\newtheorem{definition}{Definition}[section]
\newtheorem{theorem}{Theorem}[section]

\allowdisplaybreaks

\allowdisplaybreaks

\begin{document}
\title[Extended Mittag-Leffler Function and truncated {\Huge{{\Large{\Large$\nu$}}}}-fractional derivatives]
{Extended Mittag-Leffler Function and truncated {\Huge{{\Large{\Large$\nu$}}}}-fractional derivatives }

\author[A. Ghaffar, G. Rahman, K. S. Nisar,  Azeema]
{ A. Ghaffar, G. Rahman,  K. S. Nisar*,  Azeema}  % in alphabetical order
\address{A. Ghaffar\newline Department of Mathematical Science, BUITEMS Quetta, Pakistan}
\email{abdulghaffar.jaffar@gmail.com}

\address{Gauhar Rahman \newline
 Department of Mathematics, International Islamic University,
Islamabad, Pakistan}
\email{gauhar55uom@gmail.com}

\address{ K. S. Nisar\newline Department of Mathematics, College of Arts and
Science, Prince Sattam bin Abdulaziz University, Wadi Al dawaser, Riyadh
region 11991, Saudi Arabia}
\email{ksnisar1@gmail.com, n.sooppy@psau.edu.sa}

\address{Azeema \newline
Department of Mathematics, SBK Women University, Queeta, Pakistan}
\email{azeemaali9@gmail.com}
\thanks{Submitted .}

\thanks{*Corresponding author}

\begin{abstract}
The main objective of this article is to present {\Large$\nu$}-fractional derivative $\mu$-differentiable functions by considering 4-parameters extended Mittag-Leffler function (MLF). We investigate that the new {\Large$\nu$}-fractional derivative satisfies various properties of order calculus such as chain rule, product rule, Rolle's and mean-value theorems for $\mu$-differentiable function and its extension. Moreover, we define the generalized form of inverse property and the fundamental theorem of calculus and the mean-value theorem for integrals.  Also, we establish a relationship with fractional integral through truncated {\Large$\nu$}-fractional integral.
\end{abstract}

\maketitle
%%%%%%%%%%%%%%%%%%%%%%%%%%%%%%%%%
\section{Introduction}
%%%%%%%%%%%%%%%%%%%%%%%%%%%%%%%%%
During the last two decades, the interest in MLF has considerably developed. Nowadays MLF are widely used in fractional calculus with lots of interesting applications in applied sciences, engineering, special functions, probability theory,  the fractional order differential equations, and their steadily increasing importance in physics researches.

In 1903, Mittag-Leffler \cite{ Mittag} established and studies the  classic MLF with a parameter $\alpha$ and a complex variable $z$,
  which is defined as:
\begin{eqnarray}\label{1.1}
\mathbb{E}_{\alpha}(z)=\sum^{\infty}_{n=0}\frac{z^n}{\Gamma(\alpha n+1)}; \quad (\alpha \in \mathbb{R}_{0}^{+}).
\end{eqnarray}
The exponential function (EF) gives the solution of entire order differential equations (DE) with constant coefficients, but the MLF has an analogous role for solutions of no entire order DE, and can be viewed as a generalized form of the EF.

Since 1903,  numerous extensions and generalizations of the MLF have been done, when was proposed the Mittag-Leffler \cite{ Mittag}. In 1905, Wiman \cite{Wiman} introduced and discussed  2-parameters MLF.  In 1971 Prabhakar \cite{Prabhakar} proposed the so-called 3-parameters MLF, a reasonable generalization of the 2-parameters MLF. Shukla and Prajapati \cite{Shukla} in 2007 established the 4-parameters MLF. In \cite{Salim}, author introduced the 6-parameters MLF, this being a possible generalization of other MLF discussed with less than 6-parameters.

The Extended (MLF) $E^{\nu,c}_{\theta,\vartheta}(x:p)$ is recently investigated by  \cite{Ozarslan},
which is defined as:
\begin{eqnarray}\label{1.2}
\mathbb{E}^{\nu;c}_{\theta,\vartheta}(x;p)=\sum_{n=0}^{\infty}\frac{\mathbf{B}_{p}(\nu+n,c-\nu)(c)_{n}}{\mathbf{B}(\nu,c-\nu)\Gamma(\theta n+\vartheta)}\frac{x^{n}}{n!}; \,\,\,\,p>0,\Re(c)>\Re(\nu) > 0,
\end{eqnarray}
where $B_p (z,y)$ is extended beta function defined in \cite{Mainardi}  as follows:
\begin{eqnarray}\label{1.3}
B_p (z,y)=\int_0^u u^{(z-1)}(1-u)^{(y-1)}e^{-(\frac{p}{u(1-u)})}du;\,\,\,\Re(p) > 0,\,\Re(z) > 0,\,\Re(y) > 0.
\end{eqnarray}
If $p = 0$, then $B_p (z,y)$ becomes:
\begin{eqnarray}\label{1.4}
B(z,y)=\,\int_0^u\,u^{(z-1)} (1-u)^{(y-1)} du.
\end{eqnarray}

Recently, Mittal et al. \cite{Mittal} presented an extended generalized (MLF) with five positive order-parameters $\mu,\delta,\vartheta,q,c$ as:
\begin{eqnarray}\label{1.7}
\mathbb{E}^{\vartheta,q; c}_{\mu,\delta}\,(z;p)=\,\,\sum_{n=0}^{\infty}\frac{\mathbf{B}_{p}(\vartheta+nq,c-\vartheta)(\vartheta)_{nq}}{\mathbf{B}(\vartheta,c-\vartheta)\Gamma(\mu n+\delta)}\frac{z^{n}}{n!},
\end{eqnarray}
where  $\mu,\delta,\vartheta \in \mathbb{C}$ and $\Re(\mu)\, > \,0,\,\Re(\delta) > 0,\,\Re(\vartheta) > 0,q > 0\,\,$ .
%and $\,B_p (x,y)\,$ is an extended beta function defined in (\ref{2.7}).

The  extended MLF with pathway integral operator are defined by Rahman et al. (see \cite{Rahman1}).
%and find out the relevant connections of some particular cases of the main results.
 Generalized integral formulas involving the extended MLF based on the Lavoie and Trottier integral formula are  established by\cite{Rahman2}.
Nisar et al. \cite{Nisar} investigated some statistical distribution regarding fractional calculus of generalized k-MLF. The truncated {\Large{\Large$\nu$}}-fractional derivative, for $\alpha$-differentiable functions, by means of the six-parameter truncated Mittag–Leffler function was introduced in \cite{Sousa}. Motivated from the recent work of J. Vanterler da C. Sousa and E. Capelas de Oliveira \cite{Sousa}, we derive some truncated {\Large$\nu$}-fractional derivative ({\Large$\nu$}-FD)  inequalities for  $\mu$-differentiable functions through the 4-parameters extended MLF. For more about truncated $\mathcal {V} $-fractional concepts one may be referred to \cite{Sousa2, Sousa3}.
The present paper is divided into five sections. In Section 2 , we present some new truncated {\Large{\Large$\nu$}}-FD. Some  new truncated {\Large{\Large$\nu$}}-fractional integral introduced in Section 3. In section 4, we establish the truncated {\Large{\Large$\nu$}}-FD and {\Large{\Large$\nu$}}-FI) of a 4-parameters MLF. Concluding remarks close the article are discussed in section 5.

%%%%%%%%%%%%%%%%%%%%%%%%%%%%%%%%%%%%%%%%%%%%%%%%%%%%%%%%%%%%%%%%%%%%%%%%%%%%%%%%%%%%%%%%%%%%%%%%%%%%
\section{Some new truncated {\Large{\Large$\nu$}}-fractional derivative ({\Large{\Large$\nu$}}-FD)}
%%%%%%%%%%%%%%%%%%%%%%%%%%%%%%%%%%%%%%%%%%%%%%%%%%%%%%%%%%%%%%%%%%%%%%%%%%%%%%%%%%%%%%%%%%%%%%%%%%%%
Here, we propose some new truncated {\Large{\Large$\nu$}}-FD using the four parameters truncated extended MLF and we define many classical results similar to the results prevail in the entire order calculus. We derive a theorem that alludes to the law of exponents and the extension
of the $n$ order truncated {\Large{\Large$\nu$}}-FD. From these established results, we observed that the new truncated {\Large{\Large$\nu$}}-FD is linear and it fallows the product rule, the quotient rule, chain rule and the composition of $\mu$-differentiable functions.

Then, we start with the definition of the four parameters truncated MLF given by,
%\begin{eqnarray*}
%\textit{E}^{\gamma,c}_{\alpha,\beta}(z;p)=\sum^{\infty}_{n=0}\frac{B_{p}(\gamma+n,c-\gamma)}{B(\gamma,c-\gamma)}\frac{(c)_{n}}{\Gamma(\alpha_{n}+\beta)}\frac{z^{n}}{n!}
%\end{eqnarray*}
%We define truncated extended MLF as
\begin{eqnarray}\label{1}
_{i}\mathbb{E}^{\gamma,c}_{\alpha,\beta}(z;p)=\sum^{\infty}_{n=0}\frac{B_{p}(\gamma+n,c-\gamma)(c)_{n}z^{n}}{B(\gamma,c-\gamma)
{\Gamma(\alpha n+\beta)}n!},
\end{eqnarray}
where $p\geq 0$ \& $\Re(c)>\Re(\gamma)>0$.

From Eq. (\ref{1}) and  $\Gamma(\beta)$, we obtained by following truncated function, introduce by
\begin{eqnarray}\label{2}
_{i}\textit{S}^{\gamma,c}_{\alpha,\beta}(z;p)=\Gamma(\beta) _{i}\mathbb{E}^{\gamma,c}_{\alpha,\beta}(z;p)=\Gamma(\beta)\sum^{\infty}_{n=0}\frac{B_{p}(\gamma+n,c-\gamma)(c)_{n}z^{n}}{B(\gamma,c-\gamma)
{\Gamma(\alpha n+\beta)}n!},
\end{eqnarray}
\begin{definition}\label{2.1}
Let $g:[0,\infty)\rightarrow \mathbb{R}.$ For $0<\mu<1$ the new truncated {\Large{\Large$\nu$}}-FD of order $\mu$, denoted by $_{i}^{\gamma}\mathcal{V}^{\beta,c}_{\alpha,\mu}$ is defined by
\begin{eqnarray}
_{i}^{\gamma}\mathcal{V}^{\beta,c}_{\alpha,\mu}\left(g(t);p\right)=\lim_{\epsilon \rightarrow 0}\frac{g\left(t _{i}\textit{S}^{\gamma,c}_{\alpha,\beta}\left(\epsilon t^{-\mu};p\right)\right)-g(t) }{\epsilon},
\end{eqnarray}
for $\forall t>0$, $_{i}\textit{S}^{\gamma,c}_{\alpha,\beta}(.)$ is a truncated function as defined in Eq.(\ref{2}) and being $\gamma, c, \alpha, \beta \in C$ and $p\geq 0$ such that  $\Re(\gamma)>0, \Re(\alpha)>0, \Re(\beta)>0 \Re(c)>0$.\\
It is noted that if $g$ is differentiable in some $(0,a)$,
 $a>0$ and
\begin{equation}\nonumber
 \lim_{t\rightarrow 0^{+}}\left(_{i}^{\gamma}\mathcal{V}^{\beta,c}_{\alpha,\mu}\left(g(t);p\right)\right),
\end{equation}
exist, then we have
\begin{equation}\nonumber
_{i}^{\gamma}\mathcal{V}^{\beta,c}_{\alpha,\mu}\left(g(0);p\right)= \lim_{t\rightarrow 0^{+}}\left(_{i}^{\gamma}\mathcal{V}^{\beta,c}_{\alpha,\mu}\left(g(t);p\right)\right).
\end{equation}
\end{definition}
\begin{theorem}
If$ \,\,the\,function \,\,g:[0, \infty)\rightarrow \mathbb{R}\, is \,\mu-differentiable\,$ for $t_{0}>0$ and $0<\mu<1$, then $g $ is continuous in $t_{0}$.
\end{theorem}
\begin{proof}
We suppose the following notion
\begin{eqnarray}\label{4}
g\left(t_{0} \Gamma(\beta) _{i}\mathbb{E}^{\gamma,c}_{\alpha,\beta}\left(\epsilon t^{-\mu}_{0};p\right)\right)-g(t_{0})=\left(\frac{g\left(t_{0} \Gamma(\beta) _{i}\mathbb{E}^{\gamma,c}_{\alpha,\beta}\left(\epsilon t^{-\mu}_{0};p\right)\right)-g(t_{0})}{\epsilon}\right)\epsilon .
\end{eqnarray}
Now taking $\lim_{\epsilon \rightarrow 0}$ on both sides of Eq. (\ref{4}), we get
\begin{align*}
\lim_{\epsilon \rightarrow 0}g\left(t_{0} \Gamma(\beta) _{i}\mathbb{E}^{\gamma,c}_{\alpha,\beta}\left(\epsilon t^{-\mu}_{0};p\right)\right)-g(t_{0})=&\lim_{\epsilon \rightarrow 0}\left(\frac{g\left(t_{0} \Gamma(\beta) _{i}\mathbb{E}^{\gamma,c}_{\alpha,\beta}\left(\epsilon t^{-\mu}_{0};p\right)\right)-g(t_{0})}{\epsilon}\right)\lim_{\epsilon \rightarrow 0}\epsilon\nonumber\\
=&_{i}^{\gamma}\mathcal{V}^{\beta,c}_{\alpha,\mu}\left(g(t);p\right)\lim_{\epsilon \rightarrow 0}\epsilon\nonumber\\
=&0.
\end{align*}
Then, $g$ is continuous in $t_{0}$.\\
The series representation of the truncated function $_{i}\textit{S}^{\gamma,c}_{\alpha,\beta}(.)$, we have,
\begin{eqnarray}\label{5}
g\left(t \Gamma(\beta) _{i}\mathbb{E}^{\gamma,c}_{\alpha,\beta}\left(\epsilon t^{-\mu}_{0};p\right)\right)=g\left[t \Gamma(\beta) \sum^{\infty}_{n=0}\frac{B_{p}(\gamma +n,c-\gamma)(c)_{n}}{B(\gamma +n,c-\gamma)\Gamma (\alpha n+\beta)}\frac{\left(\epsilon t^{-\mu}\right)^{n}}{n!}\right].
\end{eqnarray}
As $g$ is continuous and then applying  $\lim$ as $\epsilon \rightarrow 0$ on Eq. (\ref{5}), we get
\begin{align*}
\lim_{\epsilon \rightarrow 0}g\left(t \Gamma(\beta) _{i}\mathbb{E}^{\gamma,c}_{\alpha,\beta}\left(\epsilon t^{-\mu}_{0};p\right)\right)=&\lim_{\epsilon \rightarrow 0}g\left[t \Gamma(\beta) \sum^{\infty}_{n=0}\frac{B_{p}(\gamma +n,c-\gamma)(c)_{n}}{B(\gamma +n,c-\gamma)\Gamma (\alpha n+\beta)}\frac{\left(\epsilon t^{-\mu}\right)^{n}}{n!}\right]\\
=&g\left[t \Gamma(\beta)\lim_{\epsilon \rightarrow 0} \sum^{i}_{n=0}\frac{B_{p}(\gamma +n,c-\gamma)(c)_{n}}{B(\gamma +n,c-\gamma)\Gamma (\alpha n+\beta)}\frac{\left(\epsilon t^{-\mu}\right)^{n}}{n!}\right].
\end{align*}
Besides, we have
\begin{align}\label{6}
_{i}\mathbb{E}^{\gamma,c}_{\alpha,\beta}\left(\epsilon t^{-\mu}_{0};p\right)=& \sum^{i}_{n=0}\frac{B_{p}(\gamma +n,c-\gamma)(c)_{n}\left(\epsilon t^{-\mu}\right)^{n}}{B(\gamma +n,c-\gamma)\Gamma (\alpha n+\beta)n!}\nonumber\\
%=&\sum^{i}_{n=0}\frac{B_{p}(\gamma +n,c-\gamma)}{\Gamma(\gamma)\Gamma(c-\gamma)}\frac{\Gamma(c)\Gamma(c+n)}{\Gamma(c)\Gamma (\alpha_{n}+\beta)}\frac{\left(\epsilon t^{-\mu}\right)^{n}}{n!}\nonumber\\
=&\sum^{i}_{n=0}\frac{B_{p}(\gamma +n,c-\gamma)\Gamma(c+n)\left(\epsilon t^{-\mu}\right)^{n}}{\Gamma(\gamma)\Gamma(c-\gamma)\Gamma (\alpha n+\beta)n!}\nonumber\\
=&\frac{1}{\Gamma(\beta)}+\frac{B_{p}(\gamma +1,c-\gamma)}{\Gamma(\gamma)\Gamma(c-\gamma)}\frac{\Gamma(c+1)}{\Gamma (\alpha+\beta)}\frac{\left(\epsilon t^{-\mu}\right)^{1}}{n!}\nonumber\\
+&\frac{B_{p}(\gamma +1,c-\gamma)}{\Gamma(\gamma)\Gamma(c-\gamma)}\frac{\Gamma(c+2)}{\Gamma (2\alpha +\beta)}\frac{\left(\epsilon t^{-\mu}\right)^{2}}{(n-1)!}\nonumber\\
+&\frac{B_{p}(\gamma +3,c-\gamma)}{\Gamma(\gamma)\Gamma(c-\gamma)}\frac{\Gamma(c+3)}{\Gamma (3\alpha +\beta)}\frac{\left(\epsilon t^{-\mu}\right)^{3}}{(n-2)!}+...\nonumber\\
+&\frac{B_{p}(\gamma +i,c-\gamma)}{\Gamma(\gamma)\Gamma(c-\gamma)}\frac{\Gamma(c+i)}{\Gamma (i\alpha +\beta)}\frac{\left(\epsilon t^{-\mu}\right)^{i}}{(n-i-1)!}.
\end{align}
Taking $\lim$ as $ \epsilon \rightarrow 0$ on Eq. (\ref{6}), we get
\begin{equation*}
\lim_{\epsilon \rightarrow 0} \sum^{i}_{n=0}\frac{B_{p}(\gamma +n,c-\gamma)}{B(\gamma +n,c-\gamma)}\frac{(c)_{n}}{\Gamma (\alpha n+\beta)}\frac{\left(\epsilon t^{-\mu}\right)^{n}}{n!}=\frac{1}{\Gamma(\beta)}.
\end{equation*}
Then we conclude that
\begin{equation*}
g\left(t \Gamma(\beta) _{i}\textit{E}^{\gamma,c}_{\alpha,\beta}\left(\epsilon t^{-\mu};p\right)\right)=g(t).
\end{equation*}
\end{proof}
Now, we defined  the theorem that includes the main  properties of entire order
calculus. The demonstration of the chain rule, will be check by an example, which is given in next theorem. For detail, the reasoning is the similar as
described in Theorem 2 discussed in \cite{Vanterler}.
\begin{theorem}\label{thm2}
Let $0<\mu \leq 0$, $a,b \in \mathbb{R}$, $\alpha ,\beta, \gamma, c \in C$ and $p\geq 0$ such that $\Re(\alpha)>0 ,\Re(\beta)>0, \Re(\gamma)>0,\Re(c)>0$ and  $\Re(C)>0$ and $g , h $ are $\mu$-differentiable for $t>0$, then we have:
\begin{eqnarray}\label{1.1b}
_{i}^{\gamma}\mathcal{V}^{\beta,c}_{\alpha,\mu}\left(ag+bh\right)\left(t;p\right)=a\left(_{i}^{\gamma}\mathcal{V}^{\beta,c}_{\alpha,\mu}\left(g(t);p\right)\right)+b\left( _{i}^{\gamma}\mathcal{V}^{\beta,c}_{\alpha,\mu}\left(h(t);p\right)\right),
\end{eqnarray}
\begin{eqnarray}\label{1.2b}
_{i}^{\gamma}\mathcal{V}^{\beta,c}_{\alpha,\mu}\left(g.h\right)\left(t;p\right)=g(t)\left(_{i}^{\gamma}\mathcal{V}^{\beta,c}_{\alpha,\mu}
\left(h(t);p\right)\right)+h(t)\left(_{i}^{\gamma}\mathcal{V}^{\beta,c}_{\alpha,\mu}\left(g(t);p\right)\right),
\end{eqnarray}

\begin{eqnarray}\label{1.3b}
_{i}^{\gamma}\mathcal{V}^{\beta,c}_{\alpha,\mu}\left(\frac{g}{h}\right)\left(t;p\right)=\frac{h(t)\left(_{i}^{\gamma}
\mathcal{V}^{\beta,c}_{\alpha,\mu}\left(g(t);p\right)\right)-g(t)\left(_{i}^{\gamma}\mathcal{V}^{\beta,c}_{\alpha,\mu}
\left(h(t);p\right)\right)}{\left[h(t)\right]^{2}},
\end{eqnarray}
 \begin{eqnarray}\label{1.4b}
_{i}^{\gamma}\mathcal{V}^{\beta,c}_{\alpha,\mu}\left(c;p\right)=0, \text{where $g(c)=0$ is a constant}.
\end{eqnarray}

(Chain Rule) If g is differentiable, then
\begin{eqnarray}\label{1.5}
_{i}^{\gamma}\mathcal{V}^{\beta,c}_{\alpha,\mu}\left(g(t)\right)=t^{1-\mu}\Gamma(\beta) \frac{B_{p}(\gamma +1, \gamma-c)}{B(\gamma, c-\gamma)}\frac{(c)_{1}}{\Gamma(\alpha+\beta)}\frac{dg(t)}{dt},
\end{eqnarray}
being $(c)_{1}$ be the symbol of Pochhammer.
\end{theorem}

\begin{proof}
From Eq. (\ref{6}), we have
\begin{equation*}
t \Gamma(\beta) _{i}\textit{E}^{\gamma,c}_{\alpha,\beta}\left(\epsilon t^{-\mu};p\right)=t+\frac{\Gamma(\beta)}{\Gamma(\alpha+\beta)}\frac{B_{p}(\gamma+1, c-\gamma) (c)_{1}}{B(\gamma,c-\gamma)}\epsilon t^{1-\mu}+O(\epsilon^{2}),
\end{equation*}
Introducing the following change
\begin{align*}
h=&\epsilon t^{1-\mu}\left(\Gamma(\beta)\frac{B_{p}(\gamma+1, c-\gamma)}{B(\gamma,c-\gamma)}\frac{(c)_{1}}{{\Gamma(\alpha+\beta)}}+O(\epsilon)\right)\\
\Rightarrow&\epsilon=\frac{h}{t^{1-\mu}\left(\Gamma(\beta)\frac{B_{p}(\gamma+1, c-\gamma)}{B(\gamma,c-\gamma)}\frac{(c)_{1}}{{\Gamma(\alpha+\beta)}}+O(\epsilon)\right)}.
\end{align*}
We conclude that
\begin{align*}
_{i}^{\gamma}\mathcal{V}^{\beta,c}_{\alpha,\mu}\left(g(t);p\right)=&\lim_{\epsilon \rightarrow 0}\frac{\frac{g(t+h)-g(t)}{ht^{\mu-1}}}{\Gamma(\beta)\frac{B_{p}(\gamma+1, c-\gamma)}{B(\gamma,c-\gamma)}\frac{(c)_{1}}{{\Gamma(\alpha+\beta)}}\left(1+\frac{B(\gamma, c-\gamma)\Gamma(\alpha+\beta)}{\Gamma(\beta) B_{p}(\gamma+1, c-\gamma) (c)_{1}}+O(\epsilon)\right)}\\
=&\frac{t^{\mu-1}}{\frac{\Gamma(\beta) B_{p}(\gamma+1, c-\gamma) (c)_{1}}{B(\gamma, c-\gamma)\Gamma(\alpha+\beta)}}\lim_{\epsilon \rightarrow 0}\frac{\frac{g(t+h)-g(t)}{h}}{1+\frac{B(\gamma, c-\gamma)\Gamma(\alpha+\beta)}{\Gamma(\beta) B_{p}(\gamma+1, c-\gamma) (c)_{1}}+O(\epsilon)}\\
=&t^{1-\mu}\frac{\Gamma(\beta) B_{p}(\gamma+1, c-\gamma) (c)_{1}}{B(\gamma, c-\gamma)\Gamma(\alpha+\beta)}\frac{dg(t)}{dt}\\
=&t^{1-\mu}\frac{\Gamma(\beta) B_{p}(\gamma+1, c-\gamma) \Gamma(c+1)}{\Gamma(\gamma)\Gamma(c-\gamma)\Gamma(\alpha+\beta)}\frac{dg(t)}{dt},
\end{align*}
with $t>0$, $(c)_{1}=\frac{\Gamma(c+1)}{\Gamma (c)}; B(\gamma, c-\gamma)=\frac{\Gamma(\gamma)\Gamma(c-\gamma)}{\Gamma(c)}.$\\
Another property is as follows:
\begin{equation}\label{1.5b}
_{i}^{\gamma}\mathcal{V}^{\beta,c}_{\alpha,\mu}\left((g \circ  h)(t;p)\right)=g'\left(g(t)_{i}^{\gamma}\mathcal{V}^{\beta,c}_{\alpha,\mu}(h(t);p)\right),
\end{equation}
for $g$ is $\mu$-differentiable in $h(t)$.
\end{proof}

\begin{theorem}
Let $0<\mu \leq 0$, $a,b \in \mathbb{R}$, $\alpha ,\beta, \gamma, c \in C$ and $p\geq 0$ such that $\Re(\alpha)>0 ,\Re(\beta)>0, \Re(\gamma)>0,\Re(c)>0$ and  $\Re(C)>0$ and $g$ , $h$ are $\mu$-differentiable for $t>0$. Then following implication holds:

\begin{eqnarray}\label{1.6}
_{i}^{\gamma}\mathcal{V}^{\beta,c}_{\alpha,\mu}\left(e^{at};p\right)=\Gamma(\beta)\frac{B_{p}(\gamma+1, c-\gamma)}{\Gamma(\gamma)\Gamma(c-\gamma)}\frac{\Gamma(c+1)}{\Gamma(\alpha+\beta)}t^{1-\mu}a e^{at},
\end{eqnarray}
\item
\begin{eqnarray}\label{1.7b}
_{i}^{\gamma}\mathcal{V}^{\beta,c}_{\alpha,\mu}\left(\sin(at);p\right)=\Gamma(\beta)\frac{B_{p}(\gamma+1, c-\gamma)}{\Gamma(\gamma)\Gamma(c-\gamma)}\frac{\Gamma(c+1)}{\Gamma(\alpha+\beta)}t^{1-\mu}a\cos(at),
\end{eqnarray}
\item
\begin{eqnarray}\label{1.8}
_{i}^{\gamma}\mathcal{V}^{\beta,c}_{\alpha,\mu}\left(\cos(at);p\right)=\Gamma(\beta)\frac{B_{p}(\gamma+1, c-\gamma)}{\Gamma(\gamma)\Gamma(c-\gamma)}\frac{\Gamma(c+1)}{\Gamma(\alpha+\beta)}t^{1-\mu}a\sin(at),
  \end{eqnarray}
\item
\begin{eqnarray}\label{1.9}
_{i}^{\gamma}\mathcal{V}^{\beta,c}_{\alpha,\mu}\left(t^{a};p\right)=\Gamma(\beta)\frac{B_{p}(\gamma+1, c-\gamma)}{\Gamma(\gamma)\Gamma(c-\gamma)}\frac{\Gamma(c+1)}{\Gamma(\alpha+\beta)}at^{1-\mu},
\end{eqnarray}
\item
\begin{eqnarray}\label{1.10}
_{i}^{\gamma}\mathcal{V}^{\beta,c}_{\alpha,\mu}\left(\frac{t^{\alpha}}{\alpha};p\right)=\Gamma(\beta)\frac{B_{p}(\gamma+1, c-\gamma)}{\Gamma(\gamma)\Gamma(c-\gamma)}\frac{\Gamma(c-\gamma)}{\Gamma(\alpha+\beta)}.
\end{eqnarray}

\end{theorem}

\begin{theorem}
Let $0<\mu \leq 0$, $a,b \in \mathbb{R}$, $\alpha ,\beta, \gamma, c \in C$ and $\beta> 0$ such that $\Re(\alpha)>0 ,\Re(\beta)>0, \Re(\gamma)>0,\Re(c)>0$ and  $\Re(C)>0$ and g , h are $\mu$-differentiable for $t>0$. Then, following implication holds:
\begin{eqnarray}\label{1.11}
_{i}^{\gamma}\mathcal{V}^{\beta,c}_{\alpha,\mu}\left(\sin\left(\frac{t^{\mu}}{\mu}\right);p\right)=\Gamma(\beta)\frac{B_{p}(\gamma+1, c-\gamma)}{\Gamma(\gamma)\Gamma(c-\gamma)}\frac{\Gamma(c+1)}{\Gamma(\alpha+\beta)}\cos \left(\frac{t^{\mu}}{\mu}\right),
\end{eqnarray}
\begin{eqnarray}\label{1.12}
_{i}^{\gamma}\mathcal{V}^{\beta,c}_{\alpha,\mu}\left(\cos\left(\frac{t^{\mu}}{\mu}\right);p\right)=\Gamma(\beta)\frac{B_{p}(\gamma+1, c-\gamma)}{\Gamma(\gamma)\Gamma(c-\gamma)}\frac{\Gamma(c+1)}{\Gamma(\alpha+\beta)}\sin \left(\frac{t^{\mu}}{\mu}\right),
\end{eqnarray}
\begin{eqnarray}\label{1.13}
_{i}^{\gamma}\mathcal{V}^{\beta,c}_{\alpha,\mu}\left(e^{\frac{t^{\mu}}{\mu}};p\right)=\Gamma(\beta)\frac{B_{p}(\gamma+1, c-\gamma)}{\Gamma(\gamma)\Gamma(c-\gamma)}\frac{\Gamma(c+1)}{\Gamma(\alpha+\beta)}\left(e^{\frac{t^{\mu}}{\mu}}\right).
\end{eqnarray}

\end{theorem}
Theorem below proves that commutative property depends on  fractional operator.
\begin{theorem}\label{Th5}
Let $_{i}^{\gamma}\mathcal{V}^{\beta,c}_{\alpha,\mu}\left(g(t);p\right)$ and  $_{i}^{\gamma}\mathcal{V}^{\beta,c}_{\alpha,\eta}\left(g(t);p\right)$ truncated {\Large{\Large$\nu$}}-FD derivative of the order $\mu$ $(0<\mu<1)$ and $(0<\eta<1)$ respectively. So we have
\begin{align*}
_{i}^{\gamma}\mathcal{V}^{\beta,c}_{\alpha,\mu}\left(_{i}^{\gamma}\mathcal{V}^{\beta,c}_{\alpha,\eta}\left(g(t);p\right);p\right)=&\Gamma(\beta)\frac{B_{p}(\gamma+1, c-\gamma)}{\Gamma(\gamma)\Gamma(c-\gamma)}\frac{\Gamma(c+1)}{\Gamma(\alpha+\beta)}\left[(1-\mu)\left(_{i}^{\gamma}\mathcal{V}^{\beta,c}_{\alpha,\mu}\left(g(t);p\right)\right)\right.\\
&\left. +\left(_{i}^{\gamma}\mathcal{V}^{\beta,c}_{\alpha,\eta}\left(g'(t);p\right)\right)\right].
\end{align*}
\end{theorem}
\begin{proof}
In fact, using the chain rule  (\ref{1.5}) of Theorem \ref{thm2}, we have,
\begin{align}\label{7}
_{i}^{\gamma}\mathcal{V}^{\beta,c}_{\alpha,\mu}\left(_{i}^{\gamma}\mathcal{V}^{\beta,c}_{\alpha,\eta}\left(g(t);p\right);p\right)=&_{i}^{\gamma}\mathcal{V}^{\beta,c}_{\alpha,\mu}\left(\Gamma(\beta)\frac{B_{p}(\gamma+1, c-\gamma)}{\Gamma(\gamma)\Gamma(c-\gamma)}\frac{\Gamma(c+1)}{\Gamma(\alpha+\beta)}t^{1-\eta}g'(t)\right)\nonumber\\
=&t^{1-\mu}\left(\Gamma(\beta)\frac{B_{p}(\gamma+1, c-\gamma)}{\Gamma(\gamma)\Gamma(c-\gamma)}\frac{\Gamma(c+1)}{\Gamma(\alpha+\beta)}\right)^{2}\frac{d}{dt}\left(t^{1-\eta}g'(t)\right)\nonumber\\
=&\left(\Gamma(\beta)\frac{B_{p}(\gamma+1, c-\gamma)}{\Gamma(\gamma)\Gamma(c-\gamma)}\frac{\Gamma(c+1)}{\Gamma(\alpha+\beta)}\right)^{2}\left[(1-\eta)t^{1-\mu-\eta}g'(t)\right.\nonumber\\
&\left. +t^{2-\mu-\eta}g''(t)\right]
\end{align}
By Definition (\ref{2.1}), we have
\begin{equation}\label{8}
_{i}^{\gamma}\mathcal{V}^{\beta,c}_{\alpha,\mu+\eta}\left(g(t);p\right)=\Gamma(\beta)\frac{B_{p}(\gamma+1, c-\gamma)}{\Gamma(\gamma)\Gamma(c-\gamma)}\frac{\Gamma(c+1)}{\Gamma(\alpha+\beta)}\left(t^{1-\mu-\eta }g'(t)\right)
\end{equation}
So, replacing Eq.(\ref{8}) in Eq. (\ref{7}), we conclude that
\begin{align*}
_{i}^{\gamma}\mathcal{V}^{\beta,c}_{\alpha,\mu}\left(_{i}^{\gamma}\mathcal{V}^{\beta,c}_{\alpha,\eta}\left(g(t);p\right);p\right)=&\Gamma(\beta)\frac{B_{p}(\gamma+1, c-\gamma)}{\Gamma(\gamma)\Gamma(c-\gamma)}\frac{\Gamma(c+1)}{\Gamma(\alpha+\beta)}   \left[\left(\Gamma(\beta)\frac{B_{p}(\gamma+1, c-\gamma)}{\Gamma(\gamma)\Gamma(c-\gamma)}\frac{\Gamma(c+1)}{\Gamma(\alpha+\beta)}\right)\right.\\
&\left. (1-\eta)t^{1-\mu-\eta}f'(t)+\Gamma(\beta)\frac{B_{p}(\gamma+1, c-\gamma)}{\Gamma(\gamma)\Gamma(c-\gamma)}\frac{\Gamma(c+1)}{\Gamma(\alpha+\beta)} t^{2-\mu-\eta}g''(t)\right]\\
&=\Gamma(\beta)\frac{B_{p}(\gamma+1, c-\gamma)}{\Gamma(\gamma)\Gamma(c-\gamma)}\frac{\Gamma(c+1)}{\Gamma(\alpha+\beta)}\left[(1-\eta)\left( _{i}^{\gamma}\mathcal{V}^{\beta,c}_{\alpha,\mu+\eta}\left(g(t);p\right)\right)\right.\\
&\left. +t\left( _{i}^{\gamma}\mathcal{V}^{\beta,c}_{\alpha,\mu+\eta}\left(g'(t);p\right)\right)\right].
\end{align*}
From theorem (\ref{5}), follows
\begin{equation*}
_{i}^{\gamma}\mathcal{V}^{\beta,c}_{\alpha,\mu}\left(_{i}^{\gamma}\mathcal{V}^{\beta,c}_{\alpha,\eta}\left(g(t);p\right);p\right)\neq _{i}^{\gamma}\mathcal{V}^{\beta,c}_{\alpha,\mu+\eta}\left(g(t);p\right).
\end{equation*}
\end{proof}
\begin{theorem}[Rolle's Theorem for  $\mu$-differentiable function]
Let $g:[a,b]\rightarrow \mathbb{R}$ be a function with the properties and $\mu>0$
\begin{enumerate}
\item g is $\mu$-differentiable in (a,b) for some $\mu \in [a,b]$,
\item g is continuous in $[a,b]$,
\item $g(a)=g(b)$.
\end{enumerate}
Then $\exists \,c \in (a,b), \,such\, that\, \, _{i}^{\gamma}\mathcal{V}^{\beta,c}_{\alpha,\mu}\left(g(c);p\right)=0$ with $\alpha, \beta, \gamma, c \in \mathbb{C}$ and $p \geq 0$ such that $\Re(\alpha)>0$,$\Re( \beta)>0$,$\Re(\gamma)>0$,$\Re(c)>0$ .
\end{theorem}
\begin{proof}
Since $g$ is continuous  on $[a,b]$ and $g(a)=g(b)$, there exists $c\in (a,b)$ at which the function has a local extreme. Then
\begin{align*}
_{i}^{\gamma}\mathcal{V}^{\beta,c}_{\alpha,\mu}\left(g(c);p\right)=\lim_{\epsilon \rightarrow 0^{-}}\frac{g\left(c _{i}\textit{E}^{\gamma,c}_{\alpha,\beta}\left(\epsilon c^{-\mu};p\right)\right)-g(c)}{\epsilon}=\lim_{\epsilon \rightarrow 0^{+}}\frac{g\left(c _{i}\textit{E}^{\gamma,c}_{\alpha,\beta}\left(\epsilon c^{-\mu};p\right)\right)-g(c)}{\epsilon}.
\end{align*}
But the two limits have opposite signs. Hence
\begin{equation*}
_{i}^{\gamma}\mathcal{V}^{\beta,c}_{\alpha,\mu}\left(g(c);p\right)=0.
\end{equation*}
\end{proof}
\begin{theorem}[Mean-Value Theorem (MVT) for  $\mu$-differentiable function]
Let $g:[a,b]\rightarrow \mathbb{R}$ be a function with the properties and $\mu>0$
\begin{enumerate}
\item g is continuous on $[a,b]$,
\item g is $\mu$-differentiable in $(a,b)$ for some $\mu \in (0,1)$.
\end{enumerate}
 Then, $\exists$ $c \in (a,b)$, such that
\begin{equation*}
_{i}^{\gamma}\mathcal{V}^{\beta,c}_{\alpha,\mu}\left(g(c);p\right)=\frac{g(b)-g(a)}{\frac{b^{\mu}}{\mu}-\frac{a^{\mu}}{\mu}}
\end{equation*}
where $\alpha, \beta, \gamma, c \in \mathbb{C}$ such that $\Re(\alpha)>0$,$\Re( \beta)>0$,$\Re(\gamma)>0$,$\Re(c)>0$ .

\end{theorem}
\begin{proof}
Let us consider the relation
\begin{align}\label{9a}
g(x)=g(x)-g(a)-\frac{\Gamma(\gamma)\Gamma(c-\gamma)\Gamma(\alpha+\beta)}{\Gamma(\beta)B_{p}(\gamma+1, c-\gamma)\Gamma(c+1)}\left(\frac{g(b)-g(a)}{\frac{b^{\mu}}{\mu}-\frac{a^{\mu}}{\mu}}\right)\left(\frac{1}{\mu}x^{\mu}-\frac{1}{\mu}a^{\mu}\right).
\end{align}
%The function $g$ satisfies the conditions of Rolle's theorem, then there is $c \in (a,b)$ such that $_{i}^{\gamma}\mathcal{V}^{\beta,c}_{\alpha,\mu}\left(f(c);p\right)=0$.
Taking the truncated {\Large$\nu$}-FD $_{i}^{\gamma}\mathcal{V}^{\beta,c}_{\alpha,\mu}(.)$ on both sides of Eq (\ref{9a})and the relation
$$_{i}^{\gamma}\mathcal{V}^{\beta,c}_{\alpha,\mu}\left(\frac{t^{\mu}}{\mu};p\right)=\Gamma(\beta)\frac{B_{p}(\gamma+1, c-\gamma)}{\Gamma(\gamma)\Gamma(c-\gamma)}\frac{\Gamma(c+1)}{\Gamma(\alpha+\beta)}$$
and
$_{i}^{\gamma}\mathcal{V}^{\beta,c}_{\alpha,\mu}\left(c;p\right)=0$
with c is a constant. We conclude that
\begin{equation*}
_{i}^{\gamma}\mathcal{V}^{\beta,c}_{\alpha,\mu}\left(g(c);p\right)=\frac{g(b)-g(a)}{\frac{b^{\mu}}{\mu}-\frac{a^{\mu}}{\mu}}.
\end{equation*}
\end{proof}
\begin{theorem}[Extension of MVT for $\mu$-differentiable functions]
Let $g, \,h :[a,b]\rightarrow \mathbb{R}$ be a function with the properties and $\mu>0$
\begin{enumerate}
\item g,h are continuous in [a,b],
\item g,h are $\mu$-differentiable for $\mu \in (0,1)$\end{enumerate}then,
 $\exists$ $c \in (a,b)$ such that
\begin{equation*}
\frac{_{i}^{\gamma}\mathcal{V}^{\beta,c}_{\alpha,\mu}\left(g(c);p\right)}{_{i}^{\gamma}\mathcal{V}^{\beta,c}_{\alpha,\mu}\left(h(c);p\right)}=\frac{g(b)-g(a)}{h(b)-h(a)},
\end{equation*}
where $\alpha, \beta, \gamma, c \in C$ and $p\geq 0$ such that $\Re(\alpha)>0$,$\Re( \beta)>0$,$\Re(\gamma)>0$,$\Re(c)>0$ .

\end{theorem}

\begin{definition}
Let $n<\mu \leq n+1$; $n \in N$ and $g$ $n$-differentiable for $t>0$. Then the $n^{th}$ derivative of {\Large{\Large$\nu$}}-FD is defined by
\begin{equation*}
_{i}^{\gamma}\mathcal{V}^{\beta,c}_{\alpha,\mu}\left(g(t);p\right)=\lim_{\epsilon \rightarrow 0^{-}}\frac{g^{(n)}\left(t \Gamma(\beta) _{i}\textit{E}^{\gamma,c}_{\alpha,\beta}\left(\epsilon t^{n-\mu};p\right)\right)-g^{(n)}(t)}{\epsilon},
\end{equation*}
where $\alpha, \beta, \gamma, c \in C$ and $p\geq 0$ such that $\Re(\alpha)>0$,$\Re( \beta)>0$,$\Re(\gamma)>0$,$\Re(c)>0$ , if the limit exists.
\end{definition}
%%%%%%%%%%%%%%%%%%%%%%%%%%%%%%%%%%%%%%%%%%%%%%%%%%%%%%%%%%%%%%%%%%%%%%%%%%%%%%%%
\section{{\Large{\Large$\nu$}}-fractional integral ({\Large{\Large$\nu$}}-FI)}
%%%%%%%%%%%%%%%%%%%%%%%%%%%%%%%%%%%%%%%%%%%%%%%%%%%%%%%%%%%%%%%%%%%%%%%%%%%%%%%
Here, we derive the {\Large{\Large$\nu$}}-FI of a function $g$. From the given statement,
we describe a theorem that the {\Large{\Large$\nu$}}-FI is linear, the inverse property,
the fundamental theorem of calculus (FTC), the part integration theorem, and a theorem that refer
to the mean value for integrals. Also, a theorem is given which returns the sum of the
orders of two {\Large{\Large$\nu$}}-FI, semi group property. Some other findings on the {\Large{\Large$\nu$}}-FI are also discussed.

\begin{definition}
Let $g$ be a function defined on $(a, t]$ and $0 <\mu< 1$. Also, let $t\geq a$ and $a \geq0$. Then, the {\Large{\Large$\nu$}}-FI of $g$ of order $\mu$ is stated as:
\begin{equation}\label{aa}
_{a}^{\gamma}\mathcal{I}^{\beta,c}_{\alpha,\mu}\left(g(t);p\right)=\frac{\Gamma(\gamma)\Gamma(c-\gamma)\Gamma(\alpha+\beta)}{\Gamma(\beta)B_{p}(\gamma+1, c-\gamma)\Gamma(c+1)}\int^{t}_{a}\frac{g(x)}{x^{1-\mu}}dx.
\end{equation}
where $\alpha, \beta, \gamma, c \in C$ and $p\geq 0$ such that $\Re(\alpha)>0$,$\Re( \beta)>0$,$\Re(\gamma)>0$,$\Re(c)>0$ , if the limit exists.
\end{definition}

\begin{theorem}\label{9} Let $t\geq a$ and $a \geq0$. Also, let $g,h:[a,t]\rightarrow \mathbb{R}$ be continuous functions such that $ _{a}^{\gamma}\mathcal{I}^{\beta,c}_{\alpha,\mu}\left(g(t);p\right)$, $ _{a}^{\gamma}\mathcal{I}^{\beta,c}_{\alpha,\mu}\left(h(t);p\right)$ with $0<\mu<1$, then we have
\begin{eqnarray}\label{1.14}
 _{a}^{\gamma}\mathcal{I}^{\beta,c}_{\alpha,\mu}\left(g \pm h\right)(t)= _{a}^{\gamma}\mathcal{I}^{\beta,c}_{\alpha,\mu}\left(g(t);p\right) \pm  _{a}^{\gamma}\mathcal{I}^{\beta,c}_{\alpha,\mu}\left(h(t);p\right),
\end{eqnarray}
\begin{eqnarray}\label{1.15}
_{a}^{\gamma}\mathcal{I}^{\beta,c}_{\alpha,\mu}\lambda \left(g(t);p\right)=\lambda _{a}^{\gamma}\mathcal{I}^{\beta,c}_{\alpha,\mu}\left(g(t);p\right),
\end{eqnarray}
 If $t=a$, then
\begin{eqnarray}\label{1.16}
_{a}^{\gamma}\mathcal{I}^{\beta,c}_{\alpha,\mu}\left(g(a);p\right)=0,
\end{eqnarray}
If $g(x)\geq0$, then
\begin{eqnarray}\label{1.16b}
_{a}^{\gamma}\mathcal{I}^{\beta,c}_{\alpha,\mu}\left(g(t);p\right)\geq 0.
\end{eqnarray}

\end{theorem}

\begin{theorem} Let $t\geq a$ and $a \geq0$. Also, let $g:[a,t]\rightarrow \mathbb{R}$ be continuous functions such that  $ _{a}^{\gamma}\mathcal{I}^{\beta,c}_{\alpha,\mu}\left(g(t);p\right)$, $ _{a}^{\gamma}\mathcal{I}^{\beta,c}_{\alpha,\eta}\left(g(t);p\right)$ with $0<\mu<1$ and $0<\eta<1$  then we have
\begin{align*}
_{a}^{\gamma}\mathcal{I}^{\beta,c}_{\alpha,\mu}\left(_{a}^{\gamma}\mathcal{I}^{\beta,c}_{\alpha,\eta}\left(g(t);p\right);p\right)=\frac{\Gamma(\gamma)\Gamma(c-\gamma)\Gamma(\alpha+\beta)}{\Gamma(\beta)B_{p}(\gamma+1, c-\gamma)\Gamma(c+1)}\left[\frac{t^{\mu}}{\mu}\ _{a}^{\gamma}\mathcal{I}^{\beta,c}_{\alpha,\eta}\left(g(t);p\right)-\frac{1}{\mu}\  _{a}^{\gamma}\mathcal{I}^{\beta,c}_{\alpha,\mu+\eta}\left(g(t);p\right)\right].
\end{align*}
\end{theorem}
\begin{proof}
In fact, using definition (\ref{aa}), we have
\begin{align*}
_{a}^{\gamma}\mathcal{I}^{\beta,c}_{\alpha,\mu}\left(_{a}^{\gamma}\mathcal{I}^{\beta,c}_{\alpha,\eta}\left(g(t);p\right);p\right)
=&\frac{\Gamma(\gamma)\Gamma(c-\gamma)\Gamma(\alpha+\beta)}{\Gamma(\beta)B_{p}(\gamma+1, c-\gamma)\Gamma(c+1)}\int^{t}_{a}\left( _{a}^{\gamma}\mathcal{I}^{\beta,c}_{\alpha,\eta}\left(g(t);p\right)\right)x^{\mu-1}dx\\
=&\left(\frac{\Gamma(\gamma)\Gamma(c-\gamma)\Gamma(\alpha+\beta)}{\Gamma(\beta)B_{p}(\gamma+1, c-\gamma)\Gamma(c+1)}\right)^{2} \int^{t}_{a}\left(\int^{x}_{a}\frac{g(s)}{s^{1-\eta}}ds\right)x^{\mu-1}dx\\
=&\left(\frac{\Gamma(\gamma)\Gamma(c-\gamma)\Gamma(\alpha+\beta)}{\Gamma(\beta)B_{p}(\gamma+1, c-\gamma)\Gamma(c+1)}\right)^{2} \int^{t}_{a}\frac{g(s)}{s^{1-\eta}}ds \left(\frac{t^{\mu}}{\mu}-\frac{s^{\mu}}{\mu}\right) ds\\
=&\frac{\Gamma(\gamma)\Gamma(c-\gamma)\Gamma(\alpha+\beta)}{\Gamma(\beta)B_{p}(\gamma+1, c-\gamma)\Gamma(c+1)}\left[\frac{t^{\mu}}{\mu}\ _{a}^{\gamma}\mathcal{I}^{\beta,c}_{\alpha,\eta}\left(g(t);p\right)-\frac{1}{\mu}\  _{a}^{\gamma}\mathcal{I}^{\beta,c}_{\alpha,\mu+\eta}\left(g(t);p\right)\right].
\end{align*}
From Theorem (\ref{9}), we conclude that
\begin{equation*}
_{a}^{\gamma}\mathcal{I}^{\beta,c}_{\alpha,\mu}\left(_{a}^{\gamma}\mathcal{I}^{\beta,c}_{\alpha,\eta}\left(g(t);p\right);p\right) \neq \ _{a}^{\gamma}\mathcal{I}^{\beta,c}_{\alpha,\mu+\eta}\left(g(t);p\right).
\end{equation*}
\end{proof}
\begin{theorem}[Reverse]
Let $a\geq 0$, $t\geq 0$ and $0<\mu<1$. Also, let $g$ be a continuous function such that
$_{a}^{\gamma}\mathcal{V}^{\beta,c}_{\alpha,\mu}\left(g(t);p\right)$ exist. Then
\begin{equation*}
_{a}^{\gamma}\mathcal{V}^{\beta,c}_{\alpha,\mu}\left(_{a}^{\gamma}\mathcal{I}^{\beta,c}_{\alpha,\mu}\left(g(t);p\right);p\right)=g(t).
\end{equation*}
\end{theorem}
\begin{proof}
In fact, using the chain rule and definition (\ref{aa}), we have
\begin{align*}
_{a}^{\gamma}\mathcal{V}^{\beta,c}_{\alpha,\mu}\left(_{a}^{\gamma}\mathcal{I}^{\beta,c}_{\alpha,\mu}\left(g(t);p\right);p\right)=&t^{1-\mu}\frac{\Gamma(\beta)B_{p}(\gamma+1, c-\gamma)\Gamma(c+1)}{\Gamma(\gamma)\Gamma(c-\gamma)\Gamma(\alpha+\beta)}\frac{d}{dt}\left(_{a}^{\gamma}\mathcal{I}^{\beta,c}_{\alpha,\mu}\left(g(t);p\right)\right)\\
=&g(t).
\end{align*}
\end{proof}

\begin{theorem}[FTC]
Let $g:(a,b)\rightarrow \Re$ be a differentiable function and $0<\mu<1$, then $\forall t>0$, we have
\begin{equation*}
_{a}^{\gamma}\mathcal{I}^{\beta,c}_{\alpha,\mu}\left(_{a}^{\gamma}\mathcal{V}^{\beta,c}_{\alpha,\mu}\left(g(t);p\right);p\right)=g(t)-g(a),
\end{equation*}
with $\alpha, \beta, \gamma, c \in C$ and $p\geq 0$ such that $\Re(\alpha)>0$,$\Re( \beta)>0$,$\Re(\gamma)>0$,$\Re(c)>0$ .
\end{theorem}
\begin{proof}
Applying the chain rule and the FTC for order $n$ (integer) derivatives, we have
\begin{align*}
_{a}^{\gamma}\mathcal{I}^{\beta,c}_{\alpha,\mu}\left(_{a}^{\gamma}\mathcal{V}^{\beta,c}_{\alpha,\mu}\left(g(t);p\right);p\right)=&\frac{\Gamma(\gamma)\Gamma(c-\gamma)\Gamma(\alpha+\beta)}{\Gamma(\beta)B_{p}(\gamma+1, c-\gamma)\Gamma(c+1)} \int^{t}_{a}\frac{_{a}^{\gamma}\mathcal{V}^{\beta,c}_{\alpha,\mu}\left(g(t);p\right)}{x^{1-\mu}}dx\\
=&g(t)-g(a).
\end{align*}
If $f(x)=0$ then by Eq.$ (4.2)$, we have
\begin{equation*}
_{a}^{\gamma}\mathcal{I}^{\beta,c}_{\alpha,\mu}\left(_{a}^{\gamma}\mathcal{V}^{\beta,c}_{\alpha,\mu}\left(g(t);p\right);p\right)=g(t).
\end{equation*}
\end{proof}

\begin{theorem}
Let $\alpha, \beta, \gamma, c \in C$ and $p\geq 0$ such that $R(\alpha)>0$,$\Re( \beta)>0$,$\Re(\gamma)>0$,$\Re(c)>0$ and Let $g,h:[a,b]\rightarrow \Re$ be a differentiable function and $0<\mu<1$, then  we have
\begin{equation*}
\int^{b}_{a}g(x)\left(_{a}^{\gamma}\mathcal{V}^{\beta,c}_{\alpha,\mu}\left(h(t);p\right)\right)d_{\omega}x= g(x)g(x)\Big|^{b}_{a} -\int^{b}_{a}h(x)\left(_{a}^{\gamma}\mathcal{V}^{\gamma,c}_{\alpha,\mu}\left(g(x);p\right)\right)d_{\omega}x.
\end{equation*}
with
\begin{equation*}
d_{\omega}x=\frac{\Gamma(\gamma)\Gamma(c-\gamma)\Gamma(\alpha+\beta)}{\Gamma(\beta)B_{p}(\gamma+1, c-\gamma)\Gamma(c+1)}\frac{dx}{x^{1-\mu}}.
\end{equation*}
\end{theorem}
\begin{proof}
By using the definition of {\Large \Large$\nu$}-FI Eq. (4.1), applying the chain rule and the FTC, we have
\begin{align*}
\int^{b}_{a}g(x)\left(_{a}^{\gamma}\mathcal{V}^{\beta,c}_{\alpha,\mu}\left(h(t);p\right)\right)d_{\omega}x=&\frac{\Gamma(\gamma)\Gamma(c-\gamma)\Gamma(\alpha+\beta)}{\Gamma(\beta)B_{p}(\gamma+1, c-\gamma)\Gamma(c+1)}\int^{b}_{a}g(x)\left(_{a}^{\gamma}\mathcal{V}^{\beta,c}_{\alpha,\mu}\left(h(t);p\right)\right)\frac{dx}{x^{1-\mu}}\\
=&\int^{b}_{a}g(x)h'(x)dx\\
=&g(x)h(x)\Big|^{b}_{a} -\int^{b}_{a}h(x)\left(_{a}^{\gamma}\mathcal{V}^{\gamma,c}_{\alpha,\mu}\left(g(x);p\right)\right)d_{\omega}x.
\end{align*}
\end{proof}

\begin{theorem}\label{thmb}
Let $\alpha, \beta, \gamma, c \in C$ such that $\Re(\alpha)>0$,$\Re( \beta)>0$,$\Re(\gamma)>0$,$\Re(c)>0$ and  \\
$g:[a,b]\rightarrow \Re$ be a continuous function. Then for $0<\mu<1$, we have
\begin{equation}
\Big|\ _{a}^{\gamma}\mathcal{I}^{\beta,c}_{\alpha,\mu}\left(g(t);p\right)\Big| \leq \left(_{a}^{\gamma}\mathcal{I}^{\beta,c}_{\alpha,\mu}\left(\Big|g(t)\Big|;p\right)\right).
\end{equation}
\end{theorem}
\begin{proof}
From the definition of $v$-fractional of order $\mu$, we have
\begin{align*}
\Big|\ _{a}^{\gamma}\mathcal{I}^{\beta,c}_{\alpha,\mu}\left(g(t);p\right)\Big|=&\Big|\frac{\Gamma(\gamma)\Gamma(c-\gamma)\Gamma(\alpha+\beta)}{\Gamma(\beta)B_{p}(\gamma+1, c-\gamma)\Gamma(c+1)}\int^{t}_{a}\frac{g(x)}{x^{1-\mu}}dx\Big|\\
=&\Big|\frac{\Gamma(\gamma)\Gamma(c-\gamma)\Gamma(\alpha+\beta)}{\Gamma(\beta)B_{p}(\gamma+1, c-\gamma)\Gamma(c+1)}\Big| \int^{t}_{a} \Big| \frac{g(x)}{x^{1-\mu}}\Big|dx\\
=& _{a}^{\gamma}\mathcal{I}^{\beta,c}_{\alpha,\mu} \Big|g(x)\Big|.
\end{align*}
\end{proof}

%%%%%%%%%%%%%%%%%%%%%%%%%%%%%%%%%
\begin{theorem}
Let$g:[a,b]\rightarrow \Re$ be a continuous function and $\alpha, \beta, \gamma, c \in C$ and $p>0$  such that $\Re(\alpha)>0$,$\Re( \beta)>0$,$\Re(\gamma)>0$,$\Re(c)>0$ such that
\begin{equation*}
N=\sup_{t \in [a,b]} \Big|g(x) \Big|
\end{equation*}
 Then, for all $t \in [a,b]$ and $ 0<\mu <1$, we have
 \begin{equation*}
 \Big|\ _{a}^{\gamma}\mathcal{I}^{\beta,c}_{\alpha,\mu}\left(g(t);p\right)\Big| \leq \frac{\Gamma(\gamma)\Gamma(c-\gamma)\Gamma(\alpha+\beta)}{\Gamma(\beta)B_{p}(\gamma+1, c-\gamma)\Gamma(c+1)} N\left(\frac{t^{\mu}}{\mu}-\frac{a^{\mu}}{\mu}\right)
 \end{equation*}
\end{theorem}
\begin{proof}
By Theorem (\ref{thmb}), we have
\begin{align*}
\Big|\ _{a}^{\gamma}\mathcal{I}^{\beta,c}_{\alpha,\mu}\left(g(t);p\right)\Big| \leq  _{a}^{\gamma}\mathcal{I}^{\beta,c}_{\alpha,\mu}\left|g(t);p\right| =&\frac{\Gamma(\gamma)\Gamma(c-\gamma)\Gamma(\alpha+\beta)}{\Gamma(\beta)B_{p}(\gamma+1, c-\gamma)\Gamma(c+1)} \int^{t}_{a} \Big|g(x) \Big|x^{\mu-1}dx\\
<&\frac{\Gamma(\gamma)\Gamma(c-\gamma)\Gamma(\alpha+\beta)}{\Gamma(\beta)B_{p}(\gamma+1, c-\gamma)\Gamma(c+1)} N \int^{t}_{a} x^{\mu-1}dx\\
=&\frac{\Gamma(\gamma)\Gamma(c-\gamma)\Gamma(\alpha+\beta)}{\Gamma(\beta)B_{p}(\gamma+1, c-\gamma)\Gamma(c+1)} N\left(\frac{t^{\mu}}{\mu}-\frac{a^{\mu}}{\mu}\right).
\end{align*}
\end{proof}
\begin{theorem}
Let g and h functions that satisfy the following conditions:
\begin{enumerate}
\item Continuous in [a,b].
\item Limited and integrable in [a,b].
\end{enumerate}
Besides that, let $h(x)$ a negative (or no positive) function in [a,b]. Let the set $m=\inf \{g(x) : x \in [a,b] \}$ and $M=\sup \{g(x) : x \in [a,b] \}$. Then, there exists a number $\xi \in (a,b)$ such that
\begin{equation*}
\int^{b}_{a}g(x) h(x) d_{w}x=\xi \int^{b}_{a} h(x) d_{w}x.
\end{equation*}
with
\begin{equation*}
d_{w}x=\frac{\Gamma(\gamma)\Gamma(c-\gamma)\Gamma(\alpha+\beta)}{\Gamma(\beta)B_{p}(\gamma+1, c-\gamma)\Gamma(c+1)}\frac{dx}{x^{1-\mu}}
\end{equation*}
If $g$ is continues in [a,b], then $ \exists\ x_{0} \in [a,b]$, such that
\begin{equation*}
\int^{b}_{a}g(x) h(x) d_{w}x=g(x_{0}) \int^{b}_{a} h(x) d_{w}x.
\end{equation*}
\end{theorem}
\begin{proof}
Let $m=\inf g$, $M=\sup g$ and $h(x) \geq 0$ in $[a,b]$. Then, we get
\begin{equation}\label{11}
mg(x) < g(x) h(x) < M h(x)
\end{equation}
 Multiplying by $\frac{\Gamma(\gamma)\Gamma(c-\gamma)\Gamma(\alpha+\beta)}{\Gamma(\beta)B_{p}(\gamma+1, c-\gamma)\Gamma(c+1)}$ on both sides of Eq. (\ref{11}) and integrating with respect to  $x$ on $(a,b)$, we have
 \begin{align}
 m \int^{b}_{a} \frac{\Gamma(\gamma)\Gamma(c-\gamma)\Gamma(\alpha+\beta)}{\Gamma(\beta)B_{p}(\gamma+1, c-\gamma)\Gamma(c+1)}\frac{h(x)}{x^{1-\mu}}dx <& \int^{b}_{a} \frac{\Gamma(\gamma)\Gamma(c-\gamma)\Gamma(\alpha+\beta)}{\Gamma(\beta)B_{p}(\gamma+1, c-\gamma)\Gamma(c+1)} \frac{h(x)}{x^{1-\mu}}dx\nonumber \\
 <&M \int^{b}_{a}\frac{\Gamma(\gamma)\Gamma(c-\gamma)\Gamma(\alpha+\beta)}{\Gamma(\beta)B_{p}(\gamma+1, c-\gamma)\Gamma(c+1)}\frac{h(x)}{x^{1-\mu}}dx
 \end{align}
 Then $ \exists\ x_{0} \in [a,b]$, such that
\begin{equation*}
\int^{b}_{a}g(x) f(x) d_{w}x=\xi \int^{b}_{a} h(x) d_{w}x.
\end{equation*}
with
\begin{equation*}
d_{w}x=\frac{\Gamma(\gamma)\Gamma(c-\gamma)\Gamma(\alpha+\beta)}{\Gamma(\beta)B_{p}(\gamma+1, c-\gamma)\Gamma(c+1)}\frac{dx}{x^{1-\mu}}
\end{equation*}
It is observed that when $h(x)<0$, the proof is calculated as  similar way.
In addition, by the intermediate value theorem, $g$ reaches each value in the interval $[m,M]$, then to $x_{0} \,in,\ [a,b]$, $f(x_{0})=\xi$. Then, we get
\begin{equation*}
\int^{b}_{a}g(x) h(x) d_{w}x=g(x_{0}) \int^{b}_{a} h(x) d_{w}x.
\end{equation*}
If $h(x)=0$, Eq (4.4) becomes obvious and if $h(x)>0$, then Eq (4.6) implies
\begin{equation*}
m< \frac{\int^{b}_{a}g(x) h(x) d_{\omega}x}{\int^{b}_{a}h(x)d_{\omega}x}<M
\end{equation*}
Exist, a point $x_{0} \in (a,b)$ such that $m<g(x_{0})<M$, the result follows.\\
In particular, when $h(x)=1$, by theorem (11), we have the result
\begin{align*}
\int^{b}_{a}g(x) d_{\omega}x=&g(x_{0}) \frac{\Gamma(\gamma)\Gamma(c-\gamma)\Gamma(\alpha+\beta)}{\Gamma(\beta)B_{p}(\gamma+1, c-\gamma)\Gamma(c+1)}\frac{dx}{x^{1-\mu}} \int^{b}_{a}\frac{1}{x^{1-\mu}}dx\\
=&g(x_{0})\frac{\Gamma(\gamma)\Gamma(c-\gamma)\Gamma(\alpha+\beta)}{\Gamma(\beta)B_{p}(\gamma+1, c-\gamma)\Gamma(c+1)}\left(\frac{b^{\mu}}{\mu}-\frac{a^{\mu}}{\mu}\right),
\end{align*}
this implies
\begin{align}
g(x_{0})=&\frac{\Gamma(\beta)B_{p}(\gamma+1, c-\gamma)\Gamma(c+1)}{\Gamma(\gamma)\Gamma(c-\gamma)\Gamma(\alpha+\beta)}\frac{1}{\left(\frac{b^{\mu}}{\mu}-\frac{a^{\mu}}{\mu}\right)}\int^{b}_{a}g(x)\nonumber\\
=&\frac{1}{\left(\frac{a^{\mu}}{\mu}-\frac{a^{\mu}}{\mu}\right)}\int^{b}_{a}\frac{g(x)}{x^{1-\mu}}dx.
\end{align}
\end{proof}
%%%%%%%%%%%%%%%%%%%%%%%%%%%%%%%%%%%%%%%%%%%%%%%%%%%%%%%%%%%%%%%%%%%%%%%%
\section{Derivative and integral {\Large\Large$\nu$}-fractional of a MLF}
%%%%%%%%%%%%%%%%%%%%%%%%%%%%%%%%%%%%%%%%%%%%%%%%%%%%%%%%%%%%%%%%%%%%%%%%%
MLF are very useful in the theory of fractional calculus and are important for
solving  fractional differential equations.  In this section, we obtain the {\Large\Large$\nu$}-FD of 4-parameters MLF.

\begin{theorem}\label{theorem1}
Let $ ^{\gamma}_{i}\nu^{\beta,c}_{\alpha,\mu} (g(t);p)$ the truncated {\Large{\Large$\nu$}}-FD of order $\mu$  and $\mathbb{E}_{\lambda,k}$(.) the 2-parameters MLF. Then, we have
\begin{eqnarray*}
^{\gamma}_{i}\nu^{\beta,c}_{\alpha,\mu} (\mathbb{E}_{\lambda,k}(t);p)=t^{1-\mu}\frac{\Gamma(\beta)B_{p}(\gamma+1,c-\gamma)\Gamma(c+1)}{\Gamma(\gamma)\Gamma(c-\gamma)\Gamma(\alpha+\beta)}\mathbb{E}
^{2}_{\lambda,\lambda+\delta}(t),
\end{eqnarray*}
where $\mathbb{E}^{\rho}_{\lambda,\lambda+k}(.)$ is the 3-parameters MLF.
\end{theorem}
\begin{proof} In fact , using the chain rule and the 2-parameters MLF, we have
\begin{eqnarray}\label{a}
^{\gamma}_{i}\nu^{\beta,c}_{\alpha,\mu} (\mathbb{E}_{\lambda,k}(t);p)=t^{1-\mu}\frac{\Gamma(\beta)B_{p}(\gamma+1,c-\gamma)\Gamma(c+1)}
{\Gamma(\gamma)\Gamma(c-\gamma)\Gamma(\alpha+\beta)}\frac{d}{dt}\left(\sum^{\infty}_{k=0}\frac{t^{k}}{\Gamma(\lambda k+\delta)}\right)\nonumber\\
=t^{1-\mu}\frac{\Gamma(\beta)B_{p}(\gamma+1,c-\gamma)\Gamma(c+1)}
{\Gamma(\gamma)\Gamma(c-\gamma)\Gamma(\alpha+\beta)}\sum^{\infty}_{k=0}\frac{kt^{k-1}}{\Gamma(\lambda k+\delta)}.
\end{eqnarray}
Exchanging the index, $ k\rightarrow k+1$ in  Eq. (\ref{a}),we have
\begin{align*}
^{\gamma}_{i}\nu^{\beta,c}_{\alpha,\mu} (\mathbb{E}_{\lambda,k}(t);p)&=t^{1-\mu}\frac{\Gamma(\beta)B_{p}(\gamma+1,c-\gamma)\Gamma(c+1)}
{\Gamma(\gamma)\Gamma(c-\gamma)\Gamma(\alpha+\beta)}\sum^{\infty}_{k=0}\frac{(k+1)t^{k}}{\Gamma(\lambda k+\lambda+\delta)}\\
&=t^{1-\mu}\frac{\Gamma(\beta)B_{p}(\gamma+1,c-\gamma)\Gamma(c+1)}
{\Gamma(\gamma)\Gamma(c-\gamma)\Gamma(\alpha+\beta)}\mathbb{E}^{2}_{\lambda,\lambda+\delta}(t).
\end{align*}
\end{proof}
\begin{theorem}\label{theorem3}
Let $ ^{\gamma}_{i}\nu^{\beta,c}_{\alpha,\mu} (g(t);p)$ the truncated {\Large{\Large$\nu$}}-FD of order $\mu$  and $\mathbb{E}_{\lambda,k}$(.) the 2-parameters MLF. Then, we have
\begin{eqnarray}\label{b}
^{\gamma}_{i}\nu^{\beta,c}_{\alpha,\mu} (\mathbb{E}_{\lambda,k}(t);p)=t^{1-\mu}\frac{\Gamma(\beta)B_{p}(\gamma+1,c-\gamma)\Gamma(c+1)}{\Gamma(\gamma)\Gamma(c-\gamma)\Gamma(\alpha+\beta)}\Gamma(n+2)
\mathbb{E}
^{n+2}_{\lambda,\delta+\lambda(n+1)}(t).
\end{eqnarray}
\end{theorem}
%%%%%%%%%%%%%%%%%%%%%%%%%%%%%%%%%%%%%%%%%%%%%%%%%%%%%%%
\begin{proof} Let us consider the
 result of \cite{Teodoro}, which is given below\\
\begin{eqnarray}\label{c}
\frac{d^{n}}{dt^{n}}(\mathbb{E}^{\sigma,\psi}_{\lambda,\delta}(t))=(\rho)_{q,n}\mathbb{E}^{\rho+qn,q}_{\lambda,\delta+\lambda n}(t),
\end{eqnarray}
where $\mathbb{E}^{\sigma,\psi}_{\lambda,\delta}(.)$ is the 4-parameter MLF.\\In particular, for $\rho=q=1$ in Eq. (\ref{c}), we have\\
\begin{eqnarray}\label{d}
\frac{d^{n}}{dt^{n}}(\mathbb{E}^{\sigma,\psi}_{\lambda,\delta}(t))=\Gamma(n+1)\mathbb{E}^{n+1}_{\lambda,\delta+\lambda n}(t)
\end{eqnarray}
where $\mathbb{E}^{\sigma,\psi}_{\lambda,k(.)}$ is the 3-parameters MLF.\\
Applying the entire order derivative on Eq. (\ref{d}) and choosing $q=n=1$ in Eq. (\ref{c}),\\
 \begin{eqnarray*}
 \frac{d}{dt}\left(\mathbb{E}^{n+1}_{\lambda,\delta+\lambda n}(t)\right)=\frac{\Gamma(n+2)}{\Gamma(n+1)}\mathbb{E}^{n+2}_{\lambda,\delta+\lambda(n+1)}(t),\\
 \end{eqnarray*}
 We get,
 \begin{eqnarray}\label{e}
 \frac{d^{n+1}}{dt^{n+1}}\left(\mathbb{E}_{\lambda,\delta}(t)\right)=\Gamma(n+2) \mathbb{E}^{n+2}_{\lambda,\delta+\lambda(n+1)}(t).
 \end{eqnarray}
 Using the chain rule and Eq. (\ref{e}), we conclude
  \begin{eqnarray*}
  ^{\gamma}_{i}\nu^{\beta,c;n}_{\alpha,\mu}\left(\mathbb{E}_{\lambda,\delta}(t);p\right)&=&t^{n+1-\mu}
  \frac{\Gamma(\beta)B_{p}(\gamma+1,c-\gamma)\Gamma(c+1)}{\Gamma(\gamma)\Gamma(c-\gamma)\Gamma(\alpha+\beta)}\frac{d^{n+1}}{dt^{n+1}}
  \left(\mathbb{E}_{\lambda,\delta}(t)\right)\\
  &=&t^{n+1-\mu}
  \frac{\Gamma(\beta)B_{p}(\gamma+1,c-\gamma)\Gamma(c+1)}{\Gamma(\gamma)\Gamma(c-\gamma)\Gamma(\alpha+\beta)}\Gamma(n+2)
  \mathbb{E}^{n+2}_{\lambda,\delta+\lambda(n+1)}(t).\\
\end{eqnarray*}
Let us now derive the {\Large{\Large$\nu$}}-FI of the 2-parameters MLF.
\end{proof}
%%%%%%%%%%%%%%%%%%%%%%%%%%%%%%%%%%%%%%%%%%%%%%%%%%%%%%%%%%%%%%
\begin{theorem}\label{theorem5}
Let $^{\gamma}_{a}I^{\beta,c}_{\alpha,\mu}(g(t);p)$ the {\Large{\Large$\nu$}}-FI of order $\alpha$ and $\mathbb{E}_{\lambda,\delta}(.)$, the 2-parameters
MLF. Then, we have
\begin{eqnarray*}\label{f}
^{\gamma}_{a}I^{\beta,c}_{\alpha,\mu}\left(\mathbb{E}_{\lambda,k}(t);p\right)=\frac{\Gamma(\gamma)\Gamma(c-\gamma)\Gamma(\alpha+\beta)}{\Gamma(\beta)
B_{p}(\gamma+1,c-\gamma)\Gamma(c+1)}
(t^{\mu}\mathbb{E}_{\mu+1,\delta+\mu+\alpha+1}(t)-a^{\mu}\mathbb{E}_{\lambda+1,\delta+\mu}(a)).\\
\end{eqnarray*}
\end{theorem}
\begin{proof}
Indeed  applying, the statement of {\Large{\Large$\nu$}}-FI, Eq. (\ref{aa}) and the FTC, we get
\begin{eqnarray}\label{fb}
^{\gamma}_{a}I^{\beta,c}_{\alpha,\mu}\left(\mathbb{E}_{\lambda,k}(t);p\right)&=&\frac{\Gamma(\gamma)\Gamma(c-\gamma)\Gamma(\alpha+\beta)}{\Gamma(\beta)
B_{p}(\gamma+1,c-\gamma)\Gamma(c+1)}
\int^{t}_{a}x^{\mu-1}\sum^{\infty}_{k=0}\frac{x^{k}}{\Gamma(\lambda k+\delta)}dx\nonumber\\
&=&\frac{\Gamma(\gamma)\Gamma(c-\gamma)\Gamma(\alpha+\beta)}{\Gamma(\beta)B_{p}(\gamma+1,c-\gamma)\Gamma(c+1)}
\sum^{\infty}_{k=0}\frac{1}{\Gamma(\lambda k+\delta)}\int^{t}_{a}x^{\mu+k-1}dx\nonumber\\
&=&\frac{\Gamma(\gamma)\Gamma(c-\gamma)\Gamma(\alpha+\beta)}{\Gamma(\beta)B_{p}(\gamma+1,c-\gamma)\Gamma(c+1)}
\sum^{\infty}_{k=0}\frac{1}{\Gamma(\lambda k+\delta)}\left(\frac{t^{k+\mu}}{k+\mu}-\frac{a^{k+\mu}}{k+\mu}\right)\nonumber\\
&=&\frac{\Gamma(\gamma)\Gamma(c-\gamma)\Gamma(\alpha+\beta)}{\Gamma(\beta)B_{p}(\gamma+1,c-\gamma)\Gamma(c+1)}\left(t^{\mu}\sum^{\infty}_{k=0}
\frac{t^{k}}{\Gamma((\lambda+1)k+\delta+\mu+1)}\nonumber\right.\\&&\left.-a^{\mu}\sum^{\infty}_{k=0}\frac{a^{k}}{\Gamma((\lambda+1)k+\delta+\mu+1)}\right)\\
&=&\frac{\Gamma(\gamma)\Gamma(c-\gamma)\Gamma(\alpha+\beta)}{\Gamma(\beta)B_{p}(\gamma+1,c-\gamma)\Gamma(c+1)}
\left(t^{\mu}\mathbb{E}_{\lambda+1,\delta+\mu+\alpha+1}(t)-a^{\mu}\mathbb{E}_{\lambda+1,\delta+\mu}(a)\right)\nonumber
\end{eqnarray}
Specially, applying the limit $a\rightarrow 0$, on Eq. (\ref{fb}) and substituting  $k=0$,
\begin{eqnarray}\label{g}
\lim_{a\rightarrow0}\left(a^{\mu}\sum^{\infty}_{k=0}\frac{a^{k}}{\Gamma((\lambda+1)k+\delta+\mu+1)}\right)&=&\lim_{a\rightarrow0}a^{\mu}
\frac{1}{\Gamma(\delta+\mu+1)}=0.
\end{eqnarray}
In this case, from Eq. (\ref{fb}) and Eq. (\ref{g}), we conclude that
\begin{eqnarray*}
^{\gamma}_{0}I^{\beta,c}_{\alpha,\mu}\left(\mathbb{E}_{\lambda,\delta}(t);p\right)&=&\frac{\Gamma(\gamma)\Gamma(c-\gamma)\Gamma(c+1)}
{\Gamma(\beta)B_{p}(\gamma+1,c-\gamma)\Gamma(c+1)}t^{\mu}\mathbb{E}_{\lambda+1,\delta+\mu+1}(t).\\
\end{eqnarray*}
\end{proof}
%%%%%%%%%%%%%%%%%%%%%%%%%%%%%%%%%%%%%%%%%%%%%%%%%%
\begin{theorem}\label{theorem6}
Let $ ^{\gamma}_{0}I^{\beta,c}_{\alpha,\mu}\left(g(t);p\right)$ the {\Large{\Large$\nu$}}-FI  of order $\mu,0<\mu<1$, with $a=0$ and the function \\
$g(t)=(t-x)^{\lambda}, t>x$ and $\lambda>-1$. Then we have
\begin{eqnarray*}\label{h}
^{\gamma}_{0}I^{\beta,c}_{\alpha,\mu}\left((t-x)^{\lambda};p\right)
=\frac{\Gamma(\gamma)\Gamma(c-\gamma)\Gamma(c+1)}
{\Gamma(\beta)B_{p}(\gamma+1,c-\gamma)\Gamma(c+1)}J^{\mu}t^{\lambda},
\end{eqnarray*}
where $J^{\mu}t^{\lambda}$ is the Riemann-Liouville FI of order $\mu$ \cite{Kilbas,Podlubny}.
\end{theorem}
\begin{proof} Indeed, the definition of the {\Large{\Large$\nu$}}-FI, we have
\begin{eqnarray}\label{hb}
^{\gamma}_{0}I^{\beta,c}_{\alpha,\mu}\left((t-x)^{\lambda};p\right)&=&\frac{\Gamma(\gamma)\Gamma(c-\gamma)\Gamma(\alpha+\beta)}
{\Gamma(\beta)B_{p}(\gamma+1,c-\gamma)\Gamma(c+1)}\int^{t}_{0}t^{\lambda}\left(1-\frac{x}{t}\right)^{\lambda}x^{\alpha-1}dx\nonumber\\&=&
\frac{\Gamma(\gamma)\Gamma(c-\gamma)\Gamma(\alpha+\beta)}
{\Gamma(\beta)B_{p}(\gamma+1,c-\gamma)\Gamma(c+1)}\int^{1}_{0}\left(1-\lambda\right)^{\lambda}\lambda^{\alpha-1}d\lambda\nonumber\\&=&
\frac{\Gamma(\gamma)\Gamma(c-\gamma)\Gamma(\alpha+\beta)}
{\Gamma(\beta)B_{p}(\gamma+1,c-\gamma)\Gamma(c+1)}\frac{\Gamma(\lambda+1)\Gamma(\mu)}{\Gamma(\lambda+1+\mu)}.
\end{eqnarray}
Consider the following \cite{Gorenflo},
 \begin{eqnarray}\label{i}
J^{\mu}t^{\lambda}&=&\frac{\Gamma(\lambda+1)\Gamma(\mu)}{\Gamma(\lambda+1+\mu)}t^{\lambda+\mu},t>0
\end{eqnarray}
 and $\lambda>-1 $\\Thus, from Eq. (\ref{hb}) and Eq. (\ref{i}), we conclude that
\begin{eqnarray*}
^{\gamma}_{0}I^{\beta,c}_{\alpha,\mu}\left((t-x)^{\lambda};p\right)&=&
\frac{\Gamma(\mu)\Gamma(\gamma)\Gamma(c-\gamma)\Gamma(\alpha+\beta)}{\Gamma(\beta)B_{p}(\gamma+1,c-\gamma)\Gamma(c+1)}t^{\lambda+\mu},\\
\end{eqnarray*}
where $J^{\mu}(.)$ is the Riemann-Liouville FI of order $\mu$.
\end{proof}
%%%%%%%%%%%%%%%%%%%%%%%%%%%
\section{Conclusion}
%%%%%%%%%%%%%%%%%%%%%%%%%%%
 The  new truncated {\Large{\Large$\nu$}}-FD for $\mu$-differentiable functions using the 4-parameters truncated MLF are obtained. We conclude
that the truncated {\Large{\Large$\nu$}}-FD, in this sense of fractional derivatives, presented in this case, behaves extremely well in connection to the classical properties of entire order calculus. Moreover, it was possible through of truncated {\Large{\Large$\nu$}}-FD and {\Large{\Large$\nu$}}-FI, to establish the relations with the FD and FI in the Riemann-Liouville sense. In this context, with our fractional
derivative, it was possible to establish a useful confection with the FD mentioned, as seen in section 4.  As a future work, we will  establish the truncated k-MLF from the k-MLF [22, 23] and will obtain the generalization of FD. Also, the presented results are related to the function of a variable, in this sense, we can propose a  truncated {\Large{\Large$\nu$}}-FD  with $n$ real variables
%%%%%%%%%%%%%%%%%%%%%%%%%%%%%%%%%%%%%%%%%%%%%%%%%%%%%%%%%%%%%%%%%%%%%%%%

\end{document}